\documentclass[12pt,reqno]{amsart}

\usepackage{amsmath,amssymb}
\usepackage{times}
\usepackage{amsfonts}
\usepackage{amssymb}
\usepackage{latexsym}
\usepackage{amsmath}
\usepackage[mathcal]{eucal}	
\usepackage{amstext}
\usepackage{amsthm}
\usepackage{graphics}
\usepackage{color}
\allowdisplaybreaks
\numberwithin{equation}{section}
\numberwithin{figure}{section}

\theoremstyle{plain}
\newtheorem{theorem}{ Theorem}[section]
\newtheorem{proposition}[theorem]{ Proposition}
\newtheorem{lemma}[theorem]{ Lemma}

\newtheorem{example}[theorem]{ Example}
\newtheorem{remark}[theorem]{ Remark}
\newtheorem{definition}[theorem]{ Definition}

\newcommand{\C}{\mathbb C}

\newcommand{\R}{{\mathbb R}}

\begin{document}


\title[Resolvent Estimates in Variable Exponent Classes]
{Resolvent Estimates in Variable-Exponent Schatten Classes: Optimal Coefficients in the Hilbert-Schmidt Limit
}

\author{Afnan Maqbool}
\address{Abdus Salam School of Mathematical Sciences, Government College University Lahore, 54600 Punjab, Pakistan}
\email{afnanmaqbool10@gmail.com}

\author{Ali Raza$^\ast$}\thanks{$^\ast$Corresponding author.}
\address{Abdus Salam School of Mathematical Sciences, Government College University Lahore, 54600 Punjab, Pakistan}
\email{aliasadraza@gmail.com}

\keywords{Variable exponent, Schatten--von Neumann Classes, resolvent operator, extremal analytic Toeplitz matrices.}
\subjclass[2010]{47B10, 47A10, 47A55}

\begin{abstract}
We establish a resolvent estimate for operators belonging to a
variable-exponent Schatten--von Neumann class, expressed in terms
of their departure from normality. By employing the Szehr--Zarouf extremal family of analytic
Toeplitz matrices, we further derive an explicit estimate for
the departure from normality. In the Hilbert--Schmidt case,
this family also demonstrates that our estimate is
asymptotically sharp on the logarithmic scale.
Furthermore, whereas the classical Carleman estimate for
Hilbert--Schmidt operators yields the coefficients
$a_2=b_2=\frac{1}{2}$, our results reveal a more refined
dependence of the admissible coefficients in the resolvent
estimate on the underlying sequence of exponents. In
particular, suppose that $(p_j)$ is a non-decreasing sequence
of exponents converging to $2$, and that $c=\lim_{j\to\infty}(2-p_j)\log j$
exists in $[0,\infty]$. We show that the infimum $a_\pi$
of all admissible coefficients appearing in our resolvent
estimate is given by $a_\pi=\frac{1}{2}e^{-c/2}.$ This result demonstrates that the coefficients in the resolvent estimate depends
explicitly on the behaviour of the exponent
sequence. Consequently, our findings provide a precise
description of how the choice of variable exponents
influences the coefficients in the resolvent estimate.
\end{abstract}
\maketitle

\section{Introduction}
The history of resolvent estimates for linear operators on Hilbert spaces dates back to the early days of functional analysis. One of the earliest estimates of this kind appeared in a well-known paper by Carleman in 1921 \cite{Carleman}. In this paper, he established that, for any Hilbert--Schmidt operator $A$ on a Hilbert space $H$, the resolvent $(\xi I-A)^{-1}$ satisfies the following inequality:
\begin{equation}\label{ie1}
\left\|
\prod_{i=1}^{\infty}
\left(1-\frac{\mu_i}{\xi}\right)
\exp\left(\frac{\mu_i}{\xi}\right)
(\xi I-A)^{-1}
\right\|
\leqslant
\frac{1}{|\xi|}
\exp\left(
\frac{1}{2}
+\frac{\|A\|_2^2}{2|\xi|^2}
\right),
\end{equation}
where $\mu_i$ are the eigenvalues of $A$, counted with algebraic multiplicity, and $\|\cdot\|$ and $\|\cdot\|_2$ denote the operator norm and the Hilbert--Schmidt norm, respectively. Subsequently, this result was extended to operators belonging to the Schatten--von Neumann classes $S_p$, where $0<p<\infty$; see, for example, \cite{Dunford,Simon}.

Moreover, resolvent estimates of the form (\ref{ie1}) require complete knowledge of the eigenvalues of the operator. However, in many applications, the spectrum of an operator is not completely known. This difficulty was resolved by Gil' in a series of publications beginning in 1979 (see \cite{Gil} and the references therein), and subsequently by Dechevski and Persson \cite{Dechevski} in the mid-1990s, where they derived norm estimates for \((\xi I-A)^{-1}\) in terms of the distance of $\xi$ from the spectrum $\sigma(A)$ and the \(p\)-departure from normality of the operator \(A\).

Following the ideas of Gil' and of Dechevski and Persson, in 2004 Oscar F. Bandtlow \cite{1} derived resolvent bound estimates for operators belonging to the Schatten–von Neumann ideals $S_{p}$, $0< p <\infty$. The method adopted here is the same as that proposed by Gil', while the resulting upper bounds are similar to those of Dechevski and Persson. The estimated bound is
\begin{equation}
\|(\xi I-A)^{-1}\|\leqslant\frac{1}{d(\xi,\sigma(A))}\exp\left(c_p+\frac{d_p\nu_{p}^{p}}{(d(\xi,\sigma(A)))^p}\right),
\end{equation}
where $\nu_{p}$ is the $p$-departure from normality of $A$ and $c_p,d_p$ are constants depending on $p$.

Further, a sharp finite-dimensional resolvent estimates under some spectral constraints where obtained by Oleg Szehr and Rachid Zarouf in \emph{Maximum of the resolvent over matrices with given spectrum} \cite{4}. They identify explicit triangular Toeplitz model matrices that attain the extremal resolvent bounds.

This paper introduces a variable-exponent Schatten--von Neumann
class, a construction closely related to one previously studied
by M.~I. Gil' in \cite{2}. However, the results obtained in
\cite{2} focus primarily on establishing a relationship between
the series of absolute of the eigenvalues and the modular
of a compact operator belonging to a variable-exponent class.
In contrast, our work addresses a fundamentally different
problem: we establish resolvent estimates for operators
belonging to variable-exponent Schatten--von Neumann ideals.

To this end, we first show that if the sequence of exponents
is non-decreasing and converges to a finite limit $P$, then
the corresponding variable-exponent Schatten--von Neumann
class is continuously embedded in the Schatten--von Neumann
class $S_P$ with constant exponent $P$. By exploiting this
embedding and adapting the approach developed by O.~F. Bandtlow
in \cite{1} to establish resolvent bounds for compact operators
belonging to $S_P$, we obtain resolvent norm estimates for
operators in the variable-exponent class. Furthermore, by
building on the techniques developed in \cite{1,15}, we derive
an upper bound for the Hausdorff distance between the spectra
of two operators belonging to the same variable-exponent
Schatten--von Neumann class. 

We also investigate the sharpness of
the resolvent estimate by considering the Szehr--Zarouf extremal
family of analytic Toeplitz matrices introduced by O. Szehr and
R. Zarouf in \cite{4}. For this family, we obtain sufficiently
precise estimates of the variable-exponent departure from normality
to establish the optimality of the exponent appearing in the
resolvent estimate. In the Hilbert--Schmidt case, the same family
also demonstrates asymptotic sharpness on the logarithmic scale.

Finally, we investigate that the coefficients appearing in
our estimate depend explicitly
on the behaviour of the sequence of exponents.
In the classical estimate for Hilbert--Schmidt case by Carleman \cite{Carleman}, the corresponding
coefficients are $a_2=b_2=\frac{1}{2}$. However, in the
variable-exponent setting, a more delicate phenomenon occurs.
In particular, let $(p_j)$ be a non-decreasing sequence of
exponents converging to $2$, and suppose that $c=\lim_{j\to\infty}(2-p_j)\log j$
exists in $[0,\infty]$. We prove that the infimum $a_\pi$
of all admissible coefficients in our resolvent estimate
is given by
$
{a_\pi=\frac{1}{2}e^{-c/2}.}
$ This formula reveals a precise
dependence of the resolvent coefficient on the behaviour of the exponent sequence and
highlights a distinction between the classical
Hilbert--Schmidt setting and its variable-exponent
counterpart. The proof relies on two complementary estimates. First,
we establish an upper bound by separating finitely many
singular values. Second, we derive a lower bound
by employing the Szehr--Zarouf extremal family of analytic
Toeplitz matrices.

\section{Preliminaries}

Let $\mathcal B(H)$ denote the Banach algebra of all bounded linear operators on a separable Hilbert space $H$ under the uniform operator norm $\|\cdot\|$ having spectrum and resolvent set as $\sigma(A)$ and $\rho(A)$ respectively. Moreover an operator $A$ is said to be quasi nilpotent if $\sigma(A)=\{0\}.$ Let $S_\infty$ be the closed two sided ideal of compact operators in $\mathcal B(H)$. For $A\in S_\infty$, let $(\mu_j(A))_{j=1}^{\infty}$ and $(s_j(A))_{j=1}^{\infty}$ represent the sequence of eigenvalues and singular values of $A$ respectively such that $|\mu_1(A)|\geqslant|\mu_2(A)|\geqslant\cdots $ and $s_1(A)\geqslant s_2(A)\geqslant\cdots$. Also, recall that the variable exponent sequence space \cite{9}
$$l_{p(\cdot)}=\left\{(\upsilon_j)_{j=1}^{\infty} \;\; : \;\; \sum_{j=1}^{\infty}\left(\frac{\upsilon_j}{\lambda}\right)^{p(j)}<\infty, \;\; \text{for some} \; \lambda>0 \right\}$$ 
is a quasi-Banach space with the norm
$$\|(\upsilon_j)_{j=1}^{\infty}\|_{p(\cdot)}=\inf\left\{\lambda>0 \;\; : \;\; \sum_{j=1}^{\infty}\left(\frac{|\upsilon_j|}{\lambda}\right)^{p(j)}\leqslant 1 \right\},$$
where $p(\cdot): \mathbb{N}\rightarrow (0,\infty)$ is a bounded sequence. Moreover, M.~I.~Gil$'$, in \cite{2} considered a non-decreasing sequence $\pi=(p_j)_{j=1}^\infty$ with $p_j\ge1$, and defined the Nakano modular by
\begin{equation}\label{eq:gamma}
\gamma_{\pi}(A)=\sum_{j\ge1}\frac{(s_j(A))^{p_j}}{p_j},
\end{equation}
assuming that the series converges. Further, he showed the space $S_{\gamma, \pi}$ consists of all compact operators $A\in S_{\infty}$ for which $\gamma_{\pi}(tA)<\infty$ for all $t>0$ is a normed space along with following Luxemburg norm
\[
\|{A}\|_{\pi,G}=\inf\Big\{\lambda>0:\ \sum_{j=1}^\infty\frac{(s_j(A)/\lambda)^{p_j}}{p_j}\le 1\Big\}.
\] 
Moreover, he showed that $S_{\gamma, \pi}$ is a closed two sided ideal in the algebra of all bounded linear operators on $H$.
   
\section{The variable exponent Schatten--von Neumann class}
Assume that $1\le p_1\le p_2\le\cdots\le P<\infty,$ and $P=\sup_j p_j.$ For a compact operator $A\in S_{\infty}$, define the modular
\begin{equation}\label{eq:rho}
\rho_{\pi}(A):=\sum_{j\ge1} s_j(A)^{p_j}
\end{equation}
and consider $S_{\rho, \pi}$ is the collection of all compact operators $A\in S_{\infty}$ such that $\rho(A/\lambda)<\infty$ for some $\lambda>0.$ Moreover, for any $A\in S_{\rho, \pi},$ we define 
\[
\|{A}\|_{\pi}:=\inf\Big\{\lambda>0:\ \sum_{j\ge1}\Big(\frac{s_j(A)}{\lambda}\Big)^{p_j}\le 1\Big\},
\]
where $(s_j(A))_{j=1}^\infty$ represents the non-increasing sequence of singular values of $A$ and $\pi=(p_j)_{j=1}^\infty$ is a non-decreasing bounded sequence of numbers $p_j\ge1$. 

\begin{proposition}\label{p3}   
The weighted (\ref{eq:gamma}) and unweighted (\ref{eq:rho}) modulars define the same operator class, and
\begin{equation}\label{eq:equiv}
\|A\|_{\pi,G}\le\|A\|_{\pi}\le P\,\|A\|_{\pi,G}.
\end{equation}
In particular, $\|{\cdot}\|_{\pi}$ is a quasi-norm on the same Nakano operator ideal.
\end{proposition}

\begin{proof}
Let $A,B \in S_{\rho, \pi},$ put $\alpha=\|A\|_{\pi},
\beta=\|B\|_{\pi}.$
For $\alpha,\beta>0$, the singular-value inequality and convexity along with the $p_{2j-1},p_{2j}\geq p_j$ gives 
$$
\sum_{k=1}^{\infty}
\left(\frac{s_k(A+B)}{\alpha+\beta}\right)^{p_k}
\leq
2\frac{\alpha}{\alpha+\beta}
\sum_{j=1}^{\infty}
\left(\frac{s_j(A)}{\alpha}\right)^{p_j}+
2\frac{\beta}{\alpha+\beta}
\sum_{j=1}^{\infty}
\left(\frac{s_j(B)}{\beta}\right)^{p_j}.
$$
Using $p_j\geq 1$, we conclude that
\[
\|A+B\|_{\pi} \leq 2 (\|A\|_{\pi}+\|B\|_{\pi}).
\]   
Further, for every compact operator $A$, using \eqref{eq:gamma} and \eqref{eq:rho}, we get
\begin{equation}\label{eq:modsandwich}
\frac1P\rho_{\pi}(A)\le\gamma_{\pi}(A)\le\rho_{\pi}(A).
\end{equation}
For some $\lambda>0$, let $\rho_{\pi}(A/\lambda)\le1$ then $\gamma_{\pi}(A/\lambda)\le1.$ Hence $\|A\|_{\pi}\ge\|A\|_{\pi,G}.$ Now for some $\lambda>0,$ let $\gamma_{\pi}(A/\lambda)\le1$, then by following \eqref{eq:modsandwich}, we get $\rho_{\pi}(A/\lambda)\leq P.$ As $P\ge1$, thus $P\lambda$ is admissible for $\|A\|_{\pi}$, we have $\|A\|_{\pi}\le P\lambda$. Taking the infimum over admissible $\lambda$ gives $\|A\|_{\pi}\le P\,\|A\|_{\pi,G}.$ Further, the two-sided ideal estimate is direct. Since
$$s_j(BAC)\le\|{B}\|\,s_j(A)\,\|{C}\|,$$ homogeneity of the Luxemburg gauge gives
$\|{BAC}\|_{\pi}\le\|{B}\| \,\|A\|_{\pi}\,\|{C}\|.$
\end{proof}

\subsection{Schur decomposition}
One technique for obtaining the growth of the resolvent of an operator \(A\) is to treat \(A\) as a perturbation of a normal operator \(D\) with the same spectrum by a quasi-nilpotent operator \(N\). In finite dimensions, this follows from an argument of Henrici \cite{h}; in infinite dimensions, Oscar F. Bandtlow extended Henrici's idea and provide the decomposition of an operator $A\in S_{\infty}$ into its quasi nilpotent and normal parts as mentioned in Theorem 3.2 of \cite{1}. Since \(S_{\rho, \pi}\subseteq S_{\infty}\), for any \(A\in S_{\rho, \pi}\), by Theorem 3.2 of \cite{1}, there exists a decomposition \(A=D+N\) with \(D,N\in S_{\infty}\), where \(D\) is normal and \(\sigma(A)=\sigma(D)\), while \(N\) is quasi-nilpotent. Moreover, for any \(A\in S_{\rho, \pi}\), by Proposition \ref{p3} we have \(A\in S_{\gamma, \pi},\) thus by following Theorem 2.1 of \cite{2}, the eigenvalue sequence of \(A\), ordered by decreasing modulus, lies in the same Nakano sequence class whenever \(A\) lies in the operator ideal. The singular values of the normal operator \(D\) are the moduli of these eigenvalues. Thus \(D\) belongs to the same ideal. Since the ideal is linear, \(N=A-D\) also belongs to it.

\begin{definition}\label{d1}
Let $A\in S_{\rho, \pi}$ and $A=D+N$ is the Schur decomposition of $A$ which follows by Theorem 3.2 of \cite{1}, then
\[
\eta_{\pi}(A):=\inf\big\{\|{N}\|_{\pi}:\ A=D+N\text{ is a Schur decomposition}\big\}
\]
is called the $\pi$-departure from normality of $A$. Clearly, $\eta_{\pi}(A)=0$ if and only if $A$ is normal.
\end{definition}

\begin{lemma}\label{l1}   
For every $A$ in the variable-exponent class $S_{\rho, \pi}$, we have
\begin{equation}\label{eq:embed}
\|A\|_{S_P}\le\|A\|_{\pi},
\end{equation}
thus the class $S_{\rho, \pi}$ is embedded into the constant exponent class $S_{P}.$
\end{lemma}
\begin{proof}
Let $\lambda>0$ be such that $\sum_j\big(s_j(A)/\lambda\big)^{p_j}\le1.$ Every term is at most one, so $s_j(A)/\lambda\le1$. Since $p_j\le P$, thus
$$\sum_j\Big(\frac{s_j(A)}{\lambda}\Big)^{P}\le\sum_j\Big(\frac{s_j(A)}{\lambda}\Big)^{p_j}\le1.$$
Hence $\lambda$ is admissible for $\| \cdot \|_{S_P}$, so $\|A\|_{S_P}\le\lambda$. Taking the infimum over admissible $\lambda$ proves the result.
\end{proof}

We present a resolvent estimate for operators in the variable-exponent class 
\(S_{\rho,\pi}\), together with upper bounds for the Hausdorff distance 
between the spectra of two operators belonging to this class. These results 
are obtained as applications of Theorems 4.1 and 4.2 of \cite{1}, respectively.
\begin{theorem}\label{t1}
Under the assumptions above, for every $A$ in the variable-exponent Schatten--von Neumann idal class $S_{\rho, \pi}$ and every $\xi\in\rho(A)$,
\begin{equation}\label{eq:thm1}
\|{(\xi I-A)^{-1}}\|\le\frac1{d(\xi,\sigma(A))}\exp\Big(a_P\frac{\eta_{\pi}(A)^P}{d(\xi,\sigma(A))^P}+b_P\Big),
\end{equation}
where $a_{P}$ and $b_{P}$ are constants mentioned in Theorem 4.1 of \cite{1}.
\end{theorem}
\begin{proof}  
For any $A\in S_{\rho, \pi},$ and $\xi\in \rho(A),$ we discuss two possibilities: first if $A$ is normal then $\eta_{\pi}(A)=0,$ moreover for normal operator its resolvent operator is also normal and in that case norm of the resolvent operator is equal to its spectral radius thus $\|{(\xi I-A)^{-1}}\|\leq \frac1{d(\xi,\sigma(A))}\exp(b_P)$. Secondly, if $A$ is not normal, then by Schur decomposition $A=D+N,$ where $D$ is normal and $\sigma(D)=\sigma(A),$ while $N$ is quasi-nilpotent, further $D, N\in S_{\rho, \pi}.$ Moreover, by following Lemma \ref{l1}, we have $\|{N}\|_{S_P}\le\|{N}\|_{\pi}.$ Taking the infimum over all such decompositions, $\nu_P(A)\le\eta_{\pi}(A),$ where $\nu_P(A)$ denotes the usual Schatten-$P$ departure from normality. Thus by using the definition of operator norm and following Theorem 4.1 of \cite{1} with $p=P$, we get the required result.
\end{proof}

Recall that the Hausdorff distance $h(\cdot,\cdot)$ between two compact subsets $X$ and $Y$ of $\C$ is defined as
\[
h(X,Y)=\max\{\tilde d(X,Y),\tilde d(Y,X)\},
\]
where $\tilde d(X,Y)=\underset{{x\in X}}\sup d(x,Y)$.

\begin{theorem}\label{t2}
Let $A,B$ belong to $S_{\rho, \pi}$. Then
\begin{enumerate}
\item[(i)] $\displaystyle \tilde d(\sigma(A),\sigma(B))\le\frac{\eta_{\pi}(B)}{f_P\!\Big(\dfrac{\eta_{\pi}(B)}{\|{A-B}\|}\Big)}.$
\item[(ii)] $\displaystyle h(\sigma(A),\sigma(B))\le\frac{\max\{\eta_{\pi}(A),\eta_{\pi}(B)\}}{f_P\!\Big(\dfrac{\max\{\eta_{\pi}(A),\eta_{\pi}(B)\}}{\|{A-B}\|}\Big)}$,
\end{enumerate}   
where $f_P:\R_0^+\to\R_0^+$ is the inverse of the function $g_P:\R_0^+\to\R_0^+$ defined as $g_P(x)=x\exp\big(a_Px^P+b_P\big), \  x\geq 0,$ and $a_P,b_P$ are constants mentioned in Theorem \ref{t1}.
\end{theorem}
\begin{proof}  
Let $A=B$ then $\displaystyle \tilde d(\sigma(A),\sigma(A))=0,$ so required result holds trivially. Assume $A\neq B$ and $\mu\in \sigma(A),$ if $\mu\in \sigma(B),$ then there is nothing to prove. Consider $\mu\in\sigma(A)\cap\rho(B)$. Then by following similar procedure of Theorem 4.2 in \cite{1} along with Theorem \ref{t1}, we get the required results. Further, for the case when both $A$ and $B$ are normal we have $\displaystyle h(\sigma(A),\sigma(B))\le {\|{A-B}}\| \exp(b_P).$
\end{proof}

\section{Asymptotic Sharpness via  Szehr--Zarouf extremal family }
Let us define the fully degenerate analytic Toeplitz matrix family introduced
by O. Szehr and R. Zarouf in \cite{4}:
\begin{equation}\label{eq:T}
T_{\lambda,n}=(\lambda I+S_n)(I+\lambda S_n)^{-1},
\end{equation}
where \(S_n\) denotes the nilpotent lower shift on \(\mathbb{C}^n\) and
\(0<\lambda<1\). Equivalently, we have
\begin{equation}\label{eq:Tsum}
T_{\lambda,n}
=\lambda I+(1-\lambda^2)\sum_{k=1}^{n-1}(-\lambda)^{k-1}S_n^k.
\end{equation}
By Theorem 2.1 of \cite{4}, the matrix \(T_{\lambda,n}\) is a contraction,
its spectrum is \(\{\lambda\}\), and its minimal polynomial is
\((z-\lambda)^n\). Moreover, for \(\xi=0\), it attains the Kronecker
extremal value
\begin{equation}\label{eq:inv}
\|T_{\lambda,n}^{-1}\|=\lambda^{-n}.
\end{equation}
\begin{proposition}\label{p4}    
Let $Q_{\lambda,n}=T_{\lambda,n}-\lambda I$. For every decomposition $T_{\lambda,n}=D+N$ of the type used in Definition \ref{d1}, where $D$ is normal and $\sigma(D)=\sigma(T_{\lambda,n})$, one has $D=\lambda I, \ \ N=Q_{\lambda,n},$ further, $\eta_{\pi}(T_{\lambda,n})=\|{Q_{\lambda,n}}\|_{\pi}.$
\end{proposition}

\begin{proof}
The proof follows by using Theorem 2.1 of \cite{4} along with the fact that $\eta_{\pi}(T_{\lambda,n})=\|T_{\lambda,n} - \lambda I\|_{\pi}$.
\end{proof}
Now we calculate the singular-value identities for the nilpotent part ${Q_{\lambda,n}}$ of ${T_{\lambda,n}}$. In this direction, the following result provides the product and Hilbert--Schmidt identity for the singular values of the nilpotent part.

\begin{proposition}\label{p5}    
The non-zero singular values of ${Q_{\lambda,n}}$ satisfies
\begin{equation}\label{eq:prod}
\prod_{j=1}^{n-1}s_{j,n}=(1-\lambda^2)^{n-1}
\end{equation}
and    
\begin{equation}\label{eq:HS}
\sum_{j=1}^{n-1}s_{j,n}^2=\|T_{\lambda,n} - \lambda I\|_{S_2}^2=n-1-n\lambda^2+\lambda^{2n}.
\end{equation}
\end{proposition}
\begin{proof}
Since ${Q_{\lambda,n}}$ has rank $n-1,$ let $s_{1,n}\ge s_{2,n}\ge\cdots\ge s_{n-1, n}>0$ be its non-zero singular values. Further, by \eqref{eq:Tsum}, we have
\[
Q_{\lambda,n}=(1-\lambda^2)\sum_{k=1}^{n-1}(-\lambda)^{k-1}S_n^k .
\]
For $i>j,$ the Toeplitz expansion gives
\[
(Q_{\lambda,n})_{ij}=(1-\lambda^2)(-\lambda)^{i-j-1}
\]
and all other entries are zero. The kernel is $\operatorname{span}\{e_n\}$ and the range is $\operatorname{span}\{e_2,\dots,e_n\}$. Restrict $Q_{\lambda,n}$ to $\operatorname{span}\{e_1,\dots,e_{n-1}\}$ and view it as an operator onto $\operatorname{span}\{e_2,\dots,e_n\}$. In the natural orthonormal bases, the resulting $(n-1)\times(n-1)$ matrix is lower triangular and every diagonal entry is $(1-\lambda^2)$. Its singular values are exactly $s_{1,n},\dots,s_{n-1,n}$. Since the product of singular values is the modulus of the determinant, which is the product of its diagonal entries, so $\prod_{j=1}^{n-1}s_{j,n}=(1-\lambda^2)^{n-1}.$
For the Hilbert--Schmidt norm, the $k$-th strict lower diagonal contains $n-k$ entries of modulus $(1-\lambda^2)\lambda^{k-1}$. Thus
\begin{equation}\label{eq:HSsum}
\|{Q_{\lambda,n}}\|_{S_2}^2=(1-\lambda^2)^2\sum_{k=1}^{n-1}(n-k)\lambda^{2(k-1)} .
\end{equation}
Moreover,
\[
\sum_{k=1}^{n-1}(n-k)\lambda^{2(k-1)}=\frac{(n-1)-n\lambda^2+\lambda^{2n}}{(1-\lambda^2)^2}.
\]
Thus by following \eqref{eq:HSsum}, we get
\[
\|{Q_{\lambda,n}}\|_{S_2}^2=(n-1)-n\lambda^2+\lambda^{2n}.
\]
\end{proof}
\begin{remark}
By following Proposition \ref{p4} and Proposition \ref{p5}, we get 
\begin{equation}\label{eq:eta2}
\eta_2(T_{\lambda,n})=\sqrt{n-1-n\lambda^2+\lambda^{2n}} .
\end{equation}
\end{remark}
\begin{lemma}\label{l3}
For any $x=(x_1,\dots,x_{n-1})\in \mathbb{C}^{n-1}$, we have 
\begin{equation}\label{eq:comp}
\|{x}\|_{\ell^{p_n^+}}\le\|{x}\|_{p(\cdot)}\le\|{x}\|_{\ell^{p_n^-}},
\end{equation}
where $p_n^-=\min_{1\le j\le n-1}p_j$ and $p_n^+=\max_{1\le j\le n-1}p_j.$
\end{lemma}
\begin{proof}
Suppose $\mu>0,$ satisfies $\sum_j\big(|{x_j}|/\mu\big)^{p_j}\le1.$ Then $|{x_j}| / \mu\le1$ for every $j$. Since $p_j\le p_n^+$, thus
$
\Big({|{x_j}|}/\mu\Big)^{p_n^+}\le\Big({|{x_j}|}/\mu\Big)^{p_j}.$
After summation
\[ \sum_j\left(\frac{|{x_j}|}\mu\right)^{p_n^+}\le\sum_j\Big(\frac{|{x_j}|}\mu\Big)^{p_j}\le1 .
\]
Taking the infimum over admissible $\mu$ gives $\|{x}\|_{\ell^{p_n^+}}\le\|{x}\|_{p(\cdot)}$. A similar argument is used for the other inequality with $p_n^-\le p_j.$
\end{proof}
\begin{proposition}\label{p6} For the Szehr--Zarouf extremal family $T_{\lambda,n},$ we have the following estimation of the departure from normality
\[
\eta_{\pi}(T_{\lambda,n})\ge(1-\lambda^2)(n-1)^{1/p_n^+}
\quad\text{and}\quad
\eta_{\pi}(T_{\lambda,n})\le(1+\lambda)(n-1)^{1/p_n^-}.
\]
Moreover, if $H_{\lambda,n}:=\sqrt{n-1-n\lambda^2+\lambda^{2n}}$, then
\[
\eta_{\pi}(T_{\lambda,n})\ge(n-1)^{\min\{\frac1{p_n^+}-\frac12,\,0\}}H_{\lambda,n},
\quad
\eta_{\pi}(T_{\lambda,n})\le(n-1)^{\max\{\frac1{p_n^-}-\frac12,\,0\}}H_{\lambda,n}.
\]
\end{proposition}

\begin{proof}
By Proposition \ref{p4} and Lemma \ref{l3},
\begin{equation}\label{eq:i}
\eta_{\pi}(T_{\lambda,n})=\|{Q_{\lambda,n}}\|_{\pi}\ge\|{Q_{\lambda,n}}\|_{S_{p_n^+}} .
\end{equation}
By applying the arithmetic--geometric mean inequality to $s_{j,n}^{p_n^+}$ and following \eqref{eq:prod}, we get
\[
\sum_{j=1}^{n-1}s_{j,n}^{p_n^+}\ge(n-1)(1-\lambda^2)^{p_n^+},
\]
then \eqref{eq:i} gives
\[
\eta_{\pi}(T_{\lambda,n})\ge(1-\lambda^2)(n-1)^{1/p_n^+}.
\]   
For the upper bound, by Proposition \ref{p4} and Lemma \ref{l3}, we get
\begin{align}\label{4}
\eta_{\pi}(T_{\lambda,n})=\|{Q_{\lambda,n}}\|_{\pi}&\le\|{Q_{\lambda,n}}\|_{S_{p_n^-}}=\Big(\sum_{j=1}^{n-1}s_{j,n}^{p_n^-}\Big)^{1/p_n^-}.
\end{align}
Using the fact that $\|{Q_{\lambda,n}}\|$ is the maximum singular value, we have
\begin{align*}
\Big(\sum_{j=1}^{n-1}s_{j,n}^{p_n^-}\Big)^{1/p_n^-}&\le\Big[(n-1)\max_{1\le j\le n-1}s_{j,n}^{p_n^-}\Big]^{1/p_n^-}\\
&=(n-1)^{1/p_n^-}\|{Q_{\lambda,n}}\|
\end{align*}
By Theorem 2.1 of \cite{4}, $T_{\lambda,n}$ is a contraction, thus
$$\|{Q_{\lambda,n}}\|=\|{T_{\lambda,n}-\lambda I}\|\le1+\lambda.$$
Thus by following \eqref{4}, we get
\[
\eta_{\pi}(T_{\lambda,n})\le(n-1)^{1/p_n^-}(1+\lambda).
\]
Since for any $x\in\C^r$ and $0<p<\infty$, we know that
\begin{equation}\label{eq:pnorm}
r^{\min\{\frac1p-\frac12,0\}}\|{x}\|_2\le\|{x}\|_p\le r^{\max\{\frac1p-\frac12,0\}}\|{x}\|_2 .
\end{equation}
Using Lemma \ref{l3} and above relation for $n-1$ singular values, we get
\[
\eta_{\pi}(T_{\lambda,n})=\|{Q_{\lambda,n}}\|_{\pi}\le\|{Q_{\lambda,n}}\|_{S_{p_n^-}}\le(n-1)^{\max\{\frac1{p_n^-}-\frac12,0\}}\|{Q_{\lambda,n}}\|_{S_2}.
\]
Also, by using \eqref{eq:HS}, we get
\[
\eta_{\pi}(T_{\lambda,n})\le(n-1)^{\max\{\frac1{p_n^-}-\frac12,0\}}H_{\lambda,n}.
\]
Similarly, Lemma \ref{l3} and \eqref{eq:pnorm} for $n-1$ singular values gives 
\begin{align*}
\eta_{\pi}(T_{\lambda,n})=\|{Q_{\lambda,n}}\|_{\pi}&\ge\|{Q_{\lambda,n}}\|_{S_{p_n^+}}
\ge(n-1)^{\min\{\frac1{p_n^+}-\frac12,0\}}\|{Q_{\lambda,n}}\|_{S_2}.
\end{align*}
Hence
\[
\eta_{\pi}(T_{\lambda,n})\ge(n-1)^{\min\{\frac1{p_n^+}-\frac12,0\}}H_{\lambda,n}.\qedhere
\]
\end{proof}
The following result provides the logarithmic growth of the $\pi$-departure from normality of the class of the Szehr--Zarouf extremal family $T_{\lambda,n}.$
\begin{theorem}\label{t3}
Assume that $1\le p_1\le p_2\le\cdots, \;\;\;\;\;\; p_j\to P<\infty$ and fix $\lambda\in(0,1)$. Then
\[
\lim_{n\to\infty}\frac{\log\eta_{\pi}(T_{\lambda,n})}{\log n}=\frac1P .
\]
Equivalently, $\eta_{\pi}(T_{\lambda,n})=n^{1/P+o(1)}$.
\end{theorem}

\begin{proof}
Since the exponent sequence is non-decreasing, so $p_n^+=p_{n-1}$. Further, by Proposition \ref{p6}, we get
\[
\eta_{\pi}(T_{\lambda,n})\ge(1-\lambda^2)(n-1)^{1/p_{n-1}} .
\]
Taking logarithms on both sides and dividing by $\log n$,
\[
\frac{\log\eta_{\pi}(T_{\lambda,n})}{\log n}\ge\frac{\log(1-\lambda^2)}{\log n}+\frac1{p_{n-1}}\frac{\log(n-1)}{\log n}.
\]
Since $\log(1-\lambda^2)$ is fixed and $p_{n-1}\to P$, thus 
\begin{equation}\label{eq:liminf}
\liminf_{n\to\infty}\frac{\log\eta_{\pi}(T_{\lambda,n})}{\log n}\ge \frac{1}{P}.  
\end{equation}
For the upper bound, fix $\varepsilon\in(0,P)$. Choose $N$ such that $p_j\ge P-\varepsilon,\; j\ge N.$ Since all singular values of $Q_{\lambda,n}$ are bounded by $M:=1+\lambda,$ so take $\mu_n=C_\varepsilon\,n^{1/(P-\varepsilon)},$
where $C_\varepsilon>M$ is chosen so that $\big(M/C_\varepsilon\big)^{P-\varepsilon}\le\frac14.$ For all sufficiently large $n$, $M/\mu_n\le1,$ thus  
\[
\sum_{j=N}^{n-1}\Big(\frac{s_{j,n}(Q_{\lambda,n})}{\mu_n}\Big)^{p_j}\le n\Big(\frac M{\mu_n}\Big)^{P-\varepsilon}
=n\cdot\frac{M^{P-\varepsilon}}{C_\varepsilon^{P-\varepsilon}\,n}=\Big(\frac M{C_\varepsilon}\Big)^{P-\varepsilon}\le\frac14 .
\]
Hence, there are only finitely many indices $j<N$ and $\mu_n \to \infty$. Therefore, for all sufficiently large $n$, 
\[
\sum_{j=1}^{N-1}\Big(\frac{s_{j,n}(Q_{\lambda,n})}{\mu_n}\Big)^{p_j}\le\frac34 .
\]
Thus the full Luxemburg modular is at most $1$, and
\[
\eta_{\pi}(T_{\lambda,n})=\|{Q_{\lambda,n}}\|_{\pi}\le\mu_n=C_\varepsilon n^{1/(P-\varepsilon)} .
\]
Taking logarithm and dividing by $\log n$,
\[
\frac{\log\eta_{\pi}(T_{\lambda,n})}{\log n}\le\frac{\log\big(C_\varepsilon n^{1/(P-\varepsilon)}\big)}{\log n}=\frac{\log C_\varepsilon}{\log n}+\frac1{P-\varepsilon}\cdot\frac{\log n}{\log n}.
\]
Taking $\limsup$,
\begin{equation}\label{eq:limsup}
\limsup_{n\to\infty}\frac{\log\eta_{\pi}(T_{\lambda,n})}{\log n}\le\frac1{P-\varepsilon}.
\end{equation}
For any $\varepsilon\in(0,P)$, \eqref{eq:liminf} and \eqref{eq:limsup} implies
\[
\frac1P\le\liminf_{n\to\infty}\frac{\log\eta_{\pi}(T_{\lambda,n})}{\log n}\le\limsup_{n\to\infty}\frac{\log\eta_{\pi}(T_{\lambda,n})}{\log n}\le\frac1{P-\varepsilon}.
\]
Letting $\varepsilon\to0^+$ shows that the limit exists and
\[
\lim_{n\to\infty}\frac{\log\eta_{\pi}(T_{\lambda,n})}{\log n}=\frac1P.
\]
\end{proof}
We now present the result, which shows that the exponent $P$ appearing in the resolvent estimate of Theorem~\ref{t1} for the variable-exponent class $S_{\rho,\pi}$ is optimal. This result can be obtained in two ways. The first uses the Szehr--Zarouf extremal family $T_{\lambda,n}$. The second uses the elementary nilpotent shift $S_n$, whose $n-1$ nonzero singular values are all equal to one and which satisfies $\|(\lambda I - S_n)^{-1}\| \geq \lambda^{-n}$ for some $0<\lambda<1$. Here we follow the first approach, via the Szehr--Zarouf extremal family, to remain consistent with our previous work.

\begin{theorem}\label{t4} Assume $1\le p_1\le p_2\le\cdots, \;\;\;\;\; p_j\to P<\infty,$ and fix $\lambda\in(0,1)$. Put $\eta_{n}:=\eta_{\pi}(T_{\lambda,n})$ and suppose that for some $\alpha>0$ there are constants $C_0,C_1$ independent of $n$, such that
\[
\|{T_{\lambda,n}^{-1}}\|\le\frac1\lambda\exp\Big(C_0+C_1\Big(\frac{\eta_{n}}\lambda\Big)^{\alpha}\Big)
\]
for every $n$. Then $\alpha\ge P.$
\end{theorem}

\begin{proof}
Assume $\alpha<P$. Choose $\varepsilon>0,$ so that $\dfrac{\alpha}{P-\varepsilon}<1.$
The proof of Theorem \ref{t3} gives
\[
\eta_{n}\le C_\varepsilon n^{1/(P-\varepsilon)}
\]
for all sufficiently large $n$. Hence
\[
\eta_{n}^{\alpha}=O\big(n^{\alpha/(P-\varepsilon)}\big)=o(n).
\]
Taking logarithms in the assumed estimate yields
\begin{equation}\label{eq:contra1}
\log\|{T_{\lambda,n}^{-1}}\|\le\log\frac1\lambda+C_0+C_1\Big(\frac{\eta_{n}}\lambda\Big)^{\alpha}\le o(n).
\end{equation}
On the other hand, by following \eqref{eq:inv}, we get
\[
\log\|{T_{\lambda,n}^{-1}}\|=n\log\Big(\frac1\lambda\Big).
\]
This is a positive constant times $n$, so this grows linearly in $n$, a contradiction with \eqref{eq:contra1}, therefore $\alpha\ge P$.
\end{proof}

For a constant exponent $p_j=2$, the variable-exponent norm is the Hilbert--Schmidt norm. It follows from \cite{1} that the Carleman constants can be chosen as
\[
a_2=b_2=\tfrac12.
\]
For the Szehr--Zarouf extremal family $T_{\lambda,n}$, this gives the following explicit upper bound:
\[
\|{T_{\lambda,n}^{-1}}\|\le U_{\lambda,n}:=\frac1\lambda\exp\Big(\frac12+\frac{n-1-n\lambda^2+\lambda^{2n}}{2\lambda^2}\Big).
\]
Moreover, the following theorem does not assert equality in the upper bound for each $n$. It states the precise sharpness that follows from the model family: the upper bound and the exact inverse norm have the same leading exponential rate.
\begin{theorem}\label{t5}
Let $(\lambda_n)$ be a sequence in $(0,1)$ such that $\lambda_n\to1, \;  n(1-\lambda_n)\to+\infty.$ Then
\[
\frac{\log U_{\lambda_n,n}}{\log\|{T_{\lambda_n,n}^{-1}}\|}\to1 .
\]
Equivalently, the Hilbert--Schmidt resolvent estimate has the exact leading exponential rate on this Szehr--Zarouf extremal family.
\end{theorem}
\begin{proof} Let $\delta_n:=1-\lambda_n,$ since $\delta_n\to0^+$, thus 
$$
-\log\lambda_n=-\log(1-\delta_n)=\delta_n\Big[1+\frac{\delta_n}2+\frac{\delta_n^2}3+\cdots\Big]=\delta_n(1+o(1))
$$
Moreover, by following \eqref{eq:inv}, we get
\begin{equation}\label{eq:s2}
\log\|{T_{\lambda_n,n}^{-1}}\|=n(-\log\lambda_n)=n\delta_n(1+o(1)).
\end{equation}
Now
\begin{align*}
\log U_{\lambda_n,n}&=\log\frac1{\lambda_n}+\frac12+\frac{n-1-n\lambda_n^2+\lambda_n^{2n}}{2\lambda_n^2}\\
&=-\log\lambda_n+\frac12+\frac{n(1-\lambda_n^2)-1+\lambda_n^{2n}}{2\lambda_n^2}.
\end{align*}
Because $0<\lambda_n^{2n}\le1$ and $\lambda_n\to1$, all terms outside the factor containing $n(1-\lambda_n^2)$ are $O(1)+o(n\delta_n)$. Moreover,
\[
\frac{n(1-\lambda_n^2)}{2\lambda_n^2}=\frac{n(1-\lambda_n)(1+\lambda_n)}{2\lambda_n^2}=n\delta_n\frac{1+\lambda_n}{2\lambda_n^2}=n\delta_n(1+o(1)).
\]
The assumption $n\delta_n\to+\infty$ makes the bounded terms negligible. Hence
\begin{equation}\label{eq:s3}
\log U_{\lambda_n,n}=n\delta_n(1+o(1)).
\end{equation}
Comparing the two logarithmic asymptotics \eqref{eq:s2} and \eqref{eq:s3} gives the required result.
\end{proof}
    
\section{Infimum of admissible coefficients in resolvent estimate}
The general optimality of the constants in the Carleman inequality for
operators in the Schatten class \(S_2\) was investigated by
Dechevsky and Persson in \cite{B2}. In this section, we consider a
sequence of exponents \((p_j)\) converging to \(2\) and use it to
determine the greatest lower bound of the admissible coefficients
in the resolvent estimate established in Theorem~\ref{t1}. Our
approach relies on two complementary estimates. First, we establish
an upper bound by separating finitely many singular values. We then
derive a lower bound by employing the Szehr--Zarouf extremal family.
Finally, we combine these two estimates to establish the
result in Theorem~\ref{thm:1}.

\begin{definition}\label{d2}
Let $1 \le p_1 \le p_2 \le \cdots$ and $p_j \to 2$. A number $a>0$ is called \emph{admissible} if there is a finite constant $B_a$ such that, for every operator $A\in S_{\rho,\pi}$ and every $\xi\in\rho(A)$,
\begin{equation}\label{eq:resolvent}
\|(\xi I-A)^{-1}\|\le\frac1d\exp\Big(B_a+a\Big(\frac{\eta_\pi(A)}{d}\Big)^2\Big),\,\,\, d=d(\xi,\sigma(A)).
\end{equation}
We write $a_\pi$ for the infimum of the admissible coefficients. The constant $B_a$ may depend on $(p_j)$ and on $a$, but not on $A$, on its dimension, or on $\xi$.
\end{definition}

\begin{example}\label{r0}
By Theorem 4.1 of \cite{1}, with \(p=2\), the constants are
\(a_2=b_2=\tfrac12\). Moreover, Lemma~\ref{l1}, which establishes that
\(\nu_2(A)\leq\eta_\pi(A)\), implies that every
\(a\geq\tfrac12\) is admissible, consequently,
\(a_\pi\leq \tfrac12\).
\end{example}

Now, we start our discussion to estimate the upper bound, for this let $K \in S_2$, the regularized Carleman operator
$C(K) = \det_2 (I - K)\,(I - K)^{-1}$ is considered together with its entire extension in $K$. It
satisfies the classical inequality
\[
  \|{C(K)}\| \le \exp\Bigl( \frac{1 + \|{K}\|_{S_2}^2}{2} \Bigr).
\]
This is the Carleman inequality in the Hilbert--Schmidt class, recalled in particular by
Bandtlow in \cite{1}. We first deduce from it the following lemma.

\begin{lemma}\label{lem1}   
Let $N \in S_2$ be quasinilpotent. For $k \ge 1$, set $L_k^2 = \sum_{j=1}^{k} s_j(N)^2$. Then
\[
  \log \|{(I - N)^{-1}}\| \le \frac{1}{2} \sum_{j > k} s_j(N)^2
  + \sqrt{k}\, L_k + k \log(1 + L_k) + C_k,
\]
with $C_k = 1 + \sqrt{k} + \log(k + 1)$.
\end{lemma}

\begin{proof}
We choose a singular value decomposition of $N$ and denote by $F$ the sum of the first $k$ terms. Set $R = N - F$. Then $rank(F) \le k$, $\|{F}\|_{S_2} = L_k$, $\|{R}\|_{S_2}^2 = \sum_{j > k} s_j(N)^2$ and $\langle R, F \rangle_{S_2} = 0$. Define the
operator-valued entire function
\[
  H(x) = e^{-x tr(F)}\, C(R + xF), \qquad x \in \mathbb{C}.
\]
This function is a polynomial in $x$ of degree at most $k$. To justify this point even in infinite dimension, take finite-rank orthogonal projections $P_m$, strongly increasing to $I$,
whose ranges contain the initial and final spaces of $F$. Set $R_m = P_m R P_m$. In finite
dimension, we have
\[
  e^{-x tr(F)}  C(R_m + xF) = e^{ tr(R_m)} adj(I - R_m - xF).
\]
Each minor of $I - R_m - xF$ is a polynomial of degree at most $rank (F) \le k$. The regularized
Carleman operator depends continuously on its argument in Hilbert--Schmidt norm, uniformly on
compact sets in $x$. The polynomials thus obtained therefore converge locally uniformly to $H$,
which is again a polynomial of degree at most $k$. Set $q = (1 + L_k)^{-1}$. If $|x| = q$, orthogonality gives
\[
  \|{R + tF}\|_{S_2}^2 = \|{R}\|_{S_2}^2 + q^2 L_k^2 \le \|{R}\|_{S_2}^2 + 1.
\]
On the other hand, $|tr( F)| \le \sqrt{k}\, L_k$. The Carleman inequality therefore implies
\[
  \sup_{|x| = q} \|{H(x)}\| \le
  \exp\Bigl( 1 + \tfrac{1}{2} \|{R}\|_{S_2}^2 + \sqrt{k} \Bigr).
\]
Cauchy's inequality applied to the coefficients of a polynomial of degree at most $k$ gives
\[
  \|{H(1)}\| \le (k + 1)(1 + L_k)^k
  \exp\Bigl( 1 + \tfrac{1}{2} \|{R}\|_{S_2}^2 + \sqrt{k} \Bigr).
\]
Since $N$ is quasinilpotent, $det_2(I - N) = 1$, so $(I - N)^{-1} = e^{tr( F)} H(1)$. Taking
logarithms and using again $|tr( F)| \le \sqrt{k}\, L_k$, we obtain the stated inequality.
\end{proof}

\begin{proposition}\label{prop1}
Let $c=\lim_{j\to\infty}(2-p_j)\log j\in[0,+\infty]$ and let $a>\tfrac12e^{-c/2}$. Then there is a finite constant $B_a>0$ such that every quasi-nilpotent operator $N\in S_{\rho,\pi}$ satisfies
\[
\|(I-N)^{-1}\|\le\exp\big(B_a+a\|N\|_\pi^2\big).
\]
\end{proposition}

\begin{proof}
The case $N = 0$ is immediate. Assume $N \ne 0$ and set $r = \|{N}\|_{\pi}$. Then by Lemma \ref{l1} with $P = 2$, provides $\|{N}\|_{S_2} \le r$. Then we have
\[
  \frac{s_j(N)}{r} \le j^{-1/2} \qquad (j \ge 1).
\]
Indeed,
\[
  j \Bigl(\frac{s_j(N)}{r}\Bigr)^2 \le \sum_{i=1}^{j} \Bigl(\frac{s_i(N)}{r}\Bigr)^{p_i} \le 1.
\]
For $j > k$, it follows that
\[
  \Bigl(\frac{s_j(N)}{r}\Bigr)^2 \le j^{-(2 - p_j)/2} \Bigl(\frac{s_j(N)}{r}\Bigr)^{p_j}.
\]
Summing,
\[
  \sum_{j > k} s_j(N)^2 \le \alpha_k r^2, \qquad \alpha_k = \sup_{j > k} j^{-(2 - p_j)/2}.
\]
By definition of $c$, we have $\alpha_k \to e^{-c/2},$ if $c = +\infty$, this limit is zero.
Choose $k$ large enough that $\alpha_k/2 < a$. The Lemma \ref{lem1}, with $L_k \le r$, gives
\[
  \log \|{(I - N)^{-1}}\| \le \frac{\alpha_k}{2} r^2 + \sqrt{k}\, r + k \log(1 + r) + C_k.
\]
The integer $k$ is now fixed. The last three terms are therefore bounded by
$(a - \alpha_k/2) r^2 + B_a$ for a finite constant $B_a$. The proposition is proved.
\end{proof}

\begin{proposition}\label{pl}
For every $a>\tfrac12e^{-c/2}$, there exists a finite constant $B_a$ such that every operator $A\in S_{\rho, \pi}$ and for any $\xi\in \rho(A)$ following satisfies:
$$\|(\xi I-A)^{-1}\|\le\frac{1}{d(\xi,\sigma(A))}\exp\Big(B_a+a\Big(\frac{\eta_\pi(A)}{d(\xi,\sigma(A))}\Big)^2\Big).$$
\end{proposition}

\begin{proof}   
Let $A\in S_{\rho,\pi}$, $\xi\in\rho(A)$ and $d=d(\xi,\sigma(A)).$ Consider $A=D+N$ be a Schur decomposition as mentioned in Definition \ref{d1}. Put $K=(\xi I-D)^{-1}N,$ we get $\xi I-A=(\xi I-D)(I-K).$
As $D$ is normal with $\sigma(D)=\sigma(A)$, so $\|(\xi I-D)^{-1}\|=1/d$. Moreover, by following the Theorem 3.2 of Bandtlow\cite{1}, $K$ is quasi-nilpotent and we get $\|K\|_\pi\le\|N\|_\pi/d$. Further, Proposition~\ref{prop1} applied to $K$ gives
\[
\|(I-K)^{-1}\|\le\exp\big(B_a+a\|K\|_\pi^2\big)\le\exp\Big(B_a+a\left(\frac{\|N\|_\pi}d\right)^2\Big),
\]
and since $(\xi I-A)^{-1}=(I-K)^{-1}(\xi I-D)^{-1}$,
\[
\|(\xi I-A)^{-1}\|\le\frac1d\exp\Big(B_a+a\left(\frac{\|N\|_\pi}d\right)^2\Big).
\]
The right-hand side increases with $\|N\|_\pi$, taking the infimum over all Schur decompositions replaces $\|N\|_\pi$ by $\eta_\pi(A)$.
\end{proof}
Now we calculate the lower bound by using the Szehr--Zarouf extremal family of matrices. Let us start with the following lemma.

\begin{lemma}\label{lem:2}   
Assume $c < \infty$ and fix $0 < \lambda < 1$. Then
\[
  \lim_{n \to \infty} \frac{\eta_\pi(T_{\lambda,n})^2}{n} = e^{c/2}(1 - \lambda^2).
\]
\end{lemma}

\begin{proof}
The $n - 1$ nonzero singular values $s_{j,n}$ of $Q_{\lambda,n}$ are those of
$(1 - \lambda^2)(I_{n-1} + \lambda S_{n-1})^{-1}.$ Thus $1 - \lambda \le s_{j,n} \le 1 + \lambda$. Thus by following Proposition \ref{p5}, we get
\[
  \frac{1}{n} \sum_{j=1}^{n-1} s_{j,n}^2 \longrightarrow 1 - \lambda^2.
\]
For $t > 0$, write the sum defining the quantity $\|{Q_{\lambda,n}}\|_{\pi}$ at the point
$t\sqrt{n}$:
\[
  M_n(t) = \frac{1}{n} \sum_{j=1}^{n-1} \left(\frac{s_{j,n}}{t}\right)^{p_j}
  \exp\left( \frac{(2 - p_j) \log n}{2} \right).
\]
On the indices $\varepsilon n \le j < n$, for fixed $\varepsilon > 0$, we have uniformly $p_j \to 2$ and $(2 - p_j) \log n \to c$. The singular values lying in a compact interval of
$(0, \infty)$, each term is therefore uniformly close to $e^{c/2} s_{j,n}^2 / t^2$, before
division by $n$. The part $j < \varepsilon n$ is uniformly small as $\varepsilon \to 0$. To see
this, choose $0 < \beta < p_1/2$. The indices $j < n^\beta$ contribute at most a constant times
$n^{\beta - p_1/2}$, which tends to zero. For $n^\beta \le j < \varepsilon n$, the numbers
$(2 - p_j) \log n$ are uniformly bounded, since $(2 - p_j) \log j \to c$ and
$\log n / \log j \le 1/\beta$. Their total contribution is at most a constant times
$\varepsilon$. The first terms of the sum of squares of the singular values are controlled in the
same way. Letting $n \to \infty$ then $\varepsilon \to 0$, we obtain
\[
  M_n(t) \longrightarrow \frac{e^{c/2}(1 - \lambda^2)}{t^2}.
\]
The function $M_n(t)$ is continuous and strictly decreasing in $t$. It equals one exactly at
$t = \eta_\pi(T_{\lambda,n})/\sqrt{n}$. By comparing two fixed values of $t$, on either side of
$e^{c/4}\sqrt{1 - \lambda^2}$, we obtain the stated limit.
\end{proof}

\begin{proposition}\label{prop:3}
If $c < \infty$, every admissible coefficient $a$ in the definition of $a_\pi$ satisfies
$a \ge \tfrac{1}{2} e^{-c/2}$.
\end{proposition}

\begin{proof}
Apply the resolvent bound corresponding to $a$ to $T_{\lambda,n}$ at the point $\xi = 0$, where
$d = \lambda$. The norm of the inverse and the Lemma \ref{lem:2} gives
\[
  \log(1/\lambda) \le a\, \frac{e^{c/2}(1 - \lambda^2)}{\lambda^2}.
\]
Hence $a \ge e^{-c/2}\, \frac{\lambda^2 \log(1/\lambda)}{1 - \lambda^2}.$
Letting $\lambda \to 1$, we obtain $a \ge \tfrac{1}{2} e^{-c/2}$. This holds for every admissible
coefficient, hence also for their infimum $a_\pi$.
\end{proof}
Now by combining the upper and lower bound estimates we establish the following result which is the main contribution of this section.
\begin{theorem}\label{thm:1}
Assume $1 \le p_1 \le p_2 \le \cdots$ and $p_j \to 2$. If the limit
$c = \lim_{j \to \infty} (2 - p_j) \log j$ exists and is in $[0, +\infty]$, then $a_\pi = \frac{1}{2}\, e^{-c/2},$ with $e^{-\infty} = 0$.
\end{theorem}

\begin{proof}
By Proposition~\ref{pl}, $a_\pi\le\tfrac12e^{-c/2}, \  \ \text{for every }c\in[0,+\infty].$ Moreover, if $c<\infty$, then by Proposition~\ref{prop:3} we get that every admissible $a$ is at least $\tfrac12e^{-c/2}$, hence so is their infimum, consequently $a_\pi=\tfrac12e^{-c/2}$. If $c=+\infty$, we get $a_\pi\le0$, while $a_\pi\ge0$ by definition, thus $a_\pi=0.$
\end{proof}
\begin{example}
We consider the following three choices of the exponent sequence \((p_j)\).

\begin{enumerate}
    \item[(i)] Let $p_j=2-\frac{c}{\log(j+J)}, \  \ c>0,$ where \(J\) is chosen so that \(\log(1+J)\geq c\). Then $(2-p_j)\log j=c\frac{\log j}{\log(j+J)}
    \longrightarrow c \  \ \text{as } j\to\infty.$ Consequently,
    $
    a_\pi=\frac{1}{2}e^{-c/2}.
    $

    \item[(ii)] Let $p_j=2-\frac{1}{\log^2(j+J)},$ where \(J\) is chosen so that \(\log(1+J)\geq 1\). In this case, $(2-p_j)\log j
    =\frac{\log j}{\log(j+J)}
      \frac{1}{\log(j+J)}
    \longrightarrow 0 \  \ \text{as } j\to\infty.$ Thus, \(c=0\), and hence
    $
    a_\pi=\frac{1}{2}.
    $

    \item[(iii)] Finally, let $p_j=2-\frac{1}{\sqrt{\log(j+J)}},$
    where \(J\) is chosen so that \(\log(1+J)\geq 1\). We then have
    $
    (2-p_j)\log j
    =\frac{\log j}{\log(j+J)}
      \sqrt{\log(j+J)}
    \longrightarrow \infty
      \   \  \text{as } j\to\infty.
    $
    Therefore, \(c=+\infty\), and consequently
    $
    a_\pi=0.
    $ In this case, every strictly positive coefficient \(a>0\) is
    admissible with a suitable constant \(B_a\). However, this does
    not establish the admissibility of the coefficient \(a=0\).
\end{enumerate}
\end{example}

\noindent{\bf {Other Szehr--Zarouf model-operator constructions.}}
O.~Szehr and R.~Zarouf, \emph{Explicit counter examples to Sch\"affer's conjecture} \cite{6}, construct another explicit family of Toeplitz matrices with singleton spectrum. For any such matrix $A_n$ with $\sigma(A_n)=\{\lambda_n\}$, the same elementary argument as in Proposition~\ref{p4} gives $\eta_{\pi}(A_n)=\|A_n-\lambda_n I\|_{\pi}.$ No further asymptotic statement should be made without a separate estimate
of the singular values of $A_n-\lambda_n I$.
   
\noindent{\bf {The Banach-space case.}} The resolvent and spectral variation estimates obtained in this paper also have counterparts on finite-dimensional Banach spaces. In a forthcoming paper \cite{13}, we study the case of the fixed exponent $2$.
Let $X$ be a complex Banach space of dimension $n$. We introduce the quantity $\sum_{j=1}^n|\lambda_j(A)-\mu|^2,$ where the eigenvalues are counted with multiplicity.
We obtain the following resolvent estimate, where $d=\operatorname{dist}(z,\sigma(A))$:
$$
\|(zI-A)^{-1}\|_X \le \frac{\sqrt n}{d} \exp\left(\frac12+\frac{\Delta_X(A)^2}{2d^2}
\right),
$$
here
\[
\Delta_X(A)^2
=
\inf_{\mu\in\mathbb{C}}
\left\{
n\|A-\mu I\|_X^2
-
\sum_{j=1}^{n}|\lambda_j(A)-\mu|^2
\right\}.
\]
We also show that the coefficient $1/2$ in front of $\Delta_X(A)^2/d^2$ is optimal: it cannot be replaced by a smaller constant uniformly in the dimension. This follows from the explicit Toeplitz matrices constructed by Szehr and Zarouf in Explicit counterexamples to Schäffer's conjecture \cite{6}.
Furthermore, the recent construction of Houache, Szehr and Zarouf \cite{12}, which establishes the exact asymptotic constant $\sqrt e$ in Schäffer's inequality, provides examples for which our resolvent estimate is asymptotically sharp on the logarithmic scale, after removing the factor $\sqrt n$.
As in Theorem \ref{t2}, we also obtain a corresponding spectral variation estimate. The proofs and further details will be given in \cite{13}.

\noindent{\bf {Data availability and conflicts of interests.}}
The authors confirm that all data generated or analyzed during this study are included in this article. This work does not have any conflicts of interest.

\noindent\textbf{Funding:} This research received no external funding.

\noindent\textbf{Authors' contributions:} Both authors contributed equally to this work.
\section*{Acknowledgements}
We thank Oleg Szehr and Rachid Zarouf for their helpful feedback.

\end{document}